\documentclass{article}  
\usepackage{graphicx}
  \usepackage[all,2cell]{xy}          
\usepackage{amsmath,amsthm,color}
\usepackage{amssymb}

\usepackage{tikz-cd}
\usepackage{mathrsfs}

\usepackage{enumerate}
\usepackage{mathptmx}    
\usepackage{amsmath}
\usepackage{hyperref}
  
\hypersetup{colorlinks, 
breaklinks,             
linkcolor=blue,        
urlcolor=blue}         
\newtheorem{theorem}{Theorem}[section]
\newtheorem{defn}[theorem]{Definition}
\newtheorem{remark}[theorem]{Remark} 
\newtheorem{lemma}[theorem]{Lemma}
\newtheorem{example}[theorem]{Example}
\newtheorem{proposition}[theorem]{Proposition}
\newtheorem{corollary}[theorem]{Corollary}
\newcommand{\ra}{\longrightarrow}

\newcommand{\st}{\stackrel}

\newcommand{\ma}{\mathscr{A}}

\newcommand{\mm}{\mathcal{V}}
\newcommand{\md}{\mathscr{D}}

\newcommand{\mc}{\mathscr{C}}
\newcommand{\ms}{\mathscr{S}}

\newcommand{\ds}{\displaystyle}
\newcommand{\ti}{\times}
\newcommand{\ot}{\otimes}

\newcommand{\Res}{\mbox{Res}}

\newcommand{\oppair}[4]{%
\xymatrix{%
#1 \ar@<.5ex>[r]^-{#3} &{#2}\ar@<.5ex>[l]^-{#4}%
}}
\newcommand{\oppairi}[4]{%
\xymatrix@1{%
#1 \ar@<.5ex>[r]^{#3} &{#2}\ar@<.5ex>[l]^{#4}%
}}
\DeclareRobustCommand{\coprod}{\mathop{\text{\fakecoprod}}}
\newcommand{\fakecoprod}{%
  \sbox0{$\prod$}%
  \smash{\raisebox{\dimexpr.9625\depth-\dp0}{\scalebox{1}[-1]{$\prod$}}}%
  \vphantom{$\prod$}%
}
 
\begin{document}

\title{Reedy diagrams in V-model categories}


\author{Moncef Ghazel       \and
       Fethi Kadhi  
}


\date{ }

\maketitle
\begin{abstract}
We study the category of Reedy diagrams in a $\mm$-model category. Explicitly, we show that if $K$ is a small category, $\mm$ is a closed symmetric monoidal category and $\mc$ is a closed $\mm$-module, then the diagram category $\mm^K$ is a closed symmetric monoidal category and the diagram category $\mc^K$ is a closed $\mm^K$-module. We then prove that if further $K$ is a Reedy category, $\mm$ is a monoidal model category and $\mc$ is a $\mm$-model category, then with the Reedy model category structures, $\mm^K$ is a monoidal model category and $\mc^K$ is a $\mm^K$-model category provided that either the unit $1$ of $\mm$ is cofibrant or $\mm$ is cofibrantly generated.
\end{abstract}


\section{Introduction}\label{intro}

A Quillen model category is a category with additional structure obtained by abstracting the essential homotopy theoretic properties of the ordinary category of topological spaces. This additional structure allows one to do homotopy theory in the category in question and encodes the similarities between ordinary homotopy theory of topological spaces, simplicial homotopy theory, stable homotopy theory and homological algebra \cite{DHK,GJ,HB,H,Q}. A closed symmetric monoidal model category is a Quillen model category which is also a closed symmetric monoidal category in a compatible way. Under reasonable assumptions, monoids, algebras and modules over a given monoid in such a category inherit model category structures \cite{HH,SS,W}.  The homotopy category of a closed symmetric monoidal category  has a closed symmetric monoidal category structure \cite{H}. Similarly a $\mm$-model category over a closed symmetric monoidal model category $\mm$ is a model category $\mc$ which is also a closed $\mm$-module  in a compatible way. The homotopy category $Ho(\mc)$ of a $\mm$-model category $\mc$ is a closed module over the homotopy category $Ho(\mm)$ of $\mm$ \cite{H}. Closed symmetric monoidal model categories and  modules over them  owe their relevance to a number of recent striking results related to stable homotopy theory c.f., \cite{EKMM,HSS,SB}. A spectum is now identified as a module over the sphere spectrum $S$, which is a commutative monoid in a suitable symmetric monoidal category, while a ring spectrum is a monoid in the category of $S$-modules (spectra) with the smash product over $S$ as monoidal product.\\
\indent A Reedy category is a small category $K$ equipped with some additional structure which makes it possible to inductively construct diagrams of shape $K$ in a given category. The category of functors from a Reedy category to a given model category inherits a model category structure for which the weak equivalences are the objectwise ones, no further assumptions are made on the model category in question unlike the injective and projective model structures which also have these weak equivalences.\\ 
\indent Our objective in this paper is to combine the above notions. More precisely, given a small category $K$, a closed symmetric monoidal category $\mm$  and a closed $\mm$-module $\mc$, we show that the diagram category $\mm^K$ with the objectwise product is a closed symmetric monoidal category and the diagram category $\mc^K$ with the objectwise action of $\mm^K$ is a closed $\mm^K$-module. We then prove that if further $K$ is a Reedy category, $\mm$ is a monoidal model category which either has a cofibrant unit or is cofibrantly generated  and $\mc$ is a $\mm$-model category, then with the Reedy model structures, $\mm^K$ is a monoidal model category and $\mc^K$ is a $\mm^K$-model category. \\
\indent Our methods rely on Yoneda's Lemmas for  closed symmetric monoidal categories and modules over them. We therefore briefly review these notions in Section \ref{s2}.  These will play an essential role in investigating the categorical properties of the  functor categories $\mm^K$ and  $\mc^K$ and the objectwise action of the first on the second.   In Section \ref{s4}, we briefly review the Reedy model category structure and prove a few elementary results for later use. We also recall in the same section the notion of monoidal model categories, model categories over them and give some of their elementary properties. The main results of the paper are given in Sections \ref{ns4} and \ref{s5}.  Other properties and examples to which our results apply are given in Section \ref{s6}.\\
\indent \textbf{Acknowledgment.} We would like to thank the editor and the reviewer for their thoughtful ideas and constructive comments. Their suggestions substantially improved the quality of the manuscript.

\section{Yoneda's Lemmas}\label{s2}

 All categories considered in this paper are locally small. For any category $\mc$  and any two objects $x$, $y$  in $\mc$, the set of morphisms from $x$ to $y$ is always denoted by $\mc_0(x,y).$\\

Let $\mm$ be a bicomplete closed symmetric monoidal category with monoidal product $\otimes$, unit $1$, and internal hom functor denoted exponentially so that there is a natural bijection
$$\mm_{0}(m\otimes n,p)\cong\mm_{0}(m,p^n) \qquad  m,n,p\in\mm$$
Recall that a $\mm$-module consists of a category $\mc$ together with
 an action 
$$\begin{matrix}
\otimes :& \mm  \times \mc & \ra & \mc \\
 &(m, c) & \longmapsto  & m \otimes c
\end{matrix}$$
in the sense of \cite[Section 1]{JK}. See also  \cite[Definition 4.6.1]{H}.\\
A functor $F$ between $\mm$-modules is said to be a functor of $\mm$-modules if it preserves the action of $\mm$ in the sense of  \cite[Definition 4.1.7]{H}. 
A $\mm$-module $\mc $ is said to be closed if the action of $\mm$ on $\mc$ is two sided closed, i.e., if there exist two functors 
$$\begin{matrix}
           \begin{matrix}

              \mc(-,-):& \mc^{op} \times \mc &  \ra & \mm\\

               &(c, d) & \longmapsto  & \mc(c, d)
\end{matrix} & \begin{matrix}

               (-)^{-}: &\mm^{op} \times \mc &  \ra & \mc\\

                     &(m,c) & \longmapsto  & c^m

                     \end{matrix} &

\end{matrix}$$
together with natural isomorphisms
$$\mc_0 (m \ot c,d) \cong \mm_0(m,\mc(c,d)) \cong \mc_0(c,d^m)$$
 where $m\in\mm$ and $c, d \in \mc$.\\
Clearly, a closed action of $\mm$ on $\mc$ gives rise to a closed action of $\mm$ on $\mc^{op}$, to which we will refer as the dual action of $\mm$ on $\mc^{op}$.
Recall that a closed $\mm$-module is essentially a $\mm$-enriched tensored cotensored category \cite [Section 2]{JK}. Under this identification, a functor of $\mm$-modules is just a $\mm$-functor which preserves tensors.\\
Let $K$ be a small category. $\mm$ is a $\mm$-enriched, tensored cotensored category. By \cite[Sections 3.3, 2.2  and 3.7]{K}, (ordinary) functors from $K$ to $\mm$ form a $\mm$-enriched, tensored cotensored category$\mm^K$,  where the tensor and cotensor are computed pointwise. Therefore the following three functors
$$\begin{array}{crclccl}
 \ot:&\mm\times \mm^K & \to & \mm^K& & &\\
 &(m, M) & \longmapsto  & m \otimes M:&K&\to&\mm\\
      &  &              &             &k&\longmapsto& m\otimes M_k\\
\end{array}$$
$$\begin{array}{crcllll}
 \mm^K(-,-):&(\mm^K)^{op}\times \mm^K & \to & \mm \\
 &(M, N) & \longmapsto  & \mm^K(M,N)=\int_{k\in K} N_k^{M_k}\\ \end{array}$$

$$\begin{array}{crcllll}
 (-)^{(-)}:&\mm^{op}\times \mm^K & \to & \mm^K& & &\\
 &(m, M) & \longmapsto  & M^m:&K&\to&\mm\\
   &     &              &             &k&\longmapsto&M^m_k\\
\end{array}$$
define a closed $\mm$-module structure on $\mm^K$.\\
Define 
$$ \begin{array}{lll}
 h_k:&K\to & \mm\\
 &l\longmapsto  & \ds\coprod_{K_0(k,l)}1																	
\end{array}$$
\noindent The following Lemma is a special case of the enriched Yoneda's Lemma \cite [Section 2.4]{K}, \cite[Theorem 6.3.5]{B}   and \cite[Lemma 7.3.5]{R}.
\begin{lemma}\label{l2}{(Yoneda's Lemma for closed symmetric monoidal categories)}\\
The functors $Ev,\Gamma:\mm^K\times K\to\mm$ given by 
\begin{eqnarray*}
Ev(M,k)&=&M_k\\
\Gamma(M,k)&=&\mm^K (h_k,M)
\end{eqnarray*}                 
 are equivalent.
\end{lemma}
As a consequence, one has the following well known fact
\begin{lemma}\label{nl3}
The k-evaluation functor $  \mm^{K} \st {Ev_{k}}\ra\mm$, $k \in K$, is a right adjoint with left adjoint the functor $\mm \st {F_{k}}\ra \mm^{K}$ given by $F_{k}(m)= h_{k}\ot m$.
\end{lemma}

Define the codifferential of a functor $M: K\ra \mm$ to be the functor  
$$ \begin{array}{rrcllll}
CM:&K^{op}\ti K&\to & \mm^K&&&  \\
 &(k,l)&\longmapsto & h_k\otimes M_l	:& K& \to &\mm\\ 
&&&& q& \longmapsto &  h_k(q)\otimes M_l 												
\end{array}$$
The following Lemma is a special case of \cite[(3.72)]{K} . It justifies the above terminology and may be thought of as a categorical analog of the fundamental theorem of calculus.
\begin{lemma}\label{l0} (Yoneda coreduction for closed symmetric monoidal categories)\\
Let $ M\in \mm^{K}$.  $M$ is isomorphic to the coend of $CM$, i.e.,  $$M \cong \int^{k\in K} CM(k,k).$$ \\
\end{lemma}
Let $\mc$ be a bicomplete closed $\mm$-module. Let
$$ \begin{array}{cccl}
 \otimes:&\mm\times\mc&\to & \mc\\
 &(m,c)&\longmapsto  & m\otimes c																
\end{array}$$
$$ \begin{array}{cccl}
 \mc(-,-): &\mc^{op}\times\mc&\to & \mm\\
 &(c,d)&\longmapsto  & \mc(c,d)															
\end{array}$$
$$ \begin{array}{cccl}
 (-)^{(-)}:&\mm^{op}\times\mc&\to & \mc\\
 &(m,c)&\longmapsto  & c^m																
\end{array}$$ 
the functors which define the closed $\mm$-module structure on $\mc$. Let $\mc^K$ be the category of functors from $K$ to $\mc$. Define new functors
$$\begin{array}{rrcllll}
 \otimes:&\mm^K\times \mc & \to & \mc^K& & &\\
 &(M, c) & \longmapsto  & M\otimes c:&K&\to&\mc\\
  &      &              &             &k&\longmapsto& M_k\otimes c\\
\end{array}$$
$$\begin{array}{rrcllll}
	hom_r:&\mc^{op}\times \mc^K & \to & \mm^K& & &\\
 &(c, X) & \longmapsto  & hom_r(c,X):&K&\to&\mm\\
  &      &              &             &k&\longmapsto& \mc(c,X_k)\\																								
\end{array}$$
$$\begin{array}{rccl }
hom_l:& (\mm^K)^{op}\times \mc^K & \to & \mc \\
 &(M, X) & \longmapsto  & hom_l(M,X) =\int_{k\in K} X_k^{M_k}
\end{array}$$
One has the following fact which is an easy consequence of the weighted limit functor.
\begin{lemma}\label{l3} 
$(\otimes, hom_l, hom_r)$ is an adjunction of two
variables from $\mm^K\times \mc $ to $ \mc^K$ (in the sense of \cite[Definition 4.1.12]{H}).
\end{lemma} 
Let $h_k: K \to \mm$ be the functor defined above. The following fact is a special case of \cite[(3.71)]{K}.
\begin{lemma}\label{l4}{(Yoneda's Lemma for modules)}\\
The functors $Ev,\Gamma:\mc^K\times K\to\mc$ given by
\begin{eqnarray*}
Ev(X,k)&=&X_k\\
\Gamma(X,k)&=&hom_l(h_k,X)
\end{eqnarray*}      
are equivalent.
\end{lemma}

For $X\in\mc^K$ and $k\in K$, let 
$$ \begin{array}{lll}
 h^k:&K^{op}\to & \mm\\
 &l\longmapsto  & \ds\coprod_{K_0(l,k)}1																	
\end{array}$$
 Define the codifferential $CX$ and the differential $DX$ of $X$ to be the functors
$$ \begin{array}{rrcllll}
CX:&K^{op}\ti K&\to & \mc^K&&&  \\
 &(k,l)&\longmapsto & h_k\otimes X_l	:& K& \to &\mc\\ 
&&&& q& \longmapsto &  h_k(q)\otimes X_l 												
\end{array}$$
$$ \begin{array}{rrcllll}
 DX:&K^{op}\ti K&\to & \mc^K  &&&\\
 &(k,l)&\longmapsto & X^{h^k}_l:& K& \to&\mc\\ 
&&&& q& \longmapsto &  X^{h^k(q)}_l												
\end{array}$$
One has the following form of Yoneda which is again a special case of \cite[(3.71)]{K}.
\begin{lemma}\label{l5}{(Yoneda coreduction, reduction for modules)}\\ 
Let $X\in\mc^K$, then there are natural isomorphisms
\begin{equation}\label{eq2}
X\cong\int^{k\in K}CX(k,k)
\end{equation} 
\begin{equation}\label{eq3} 
X\cong\int_{k\in K}DX(k,k)
\end{equation}
\end{lemma}
The following adjoint functor theorem \cite[Theorem 4.51]{K}  is a consequence of Lemmas \ref{l4} and \ref{l5}. 
\begin{theorem}\label{nt1}
Let $\mm^{K} \st {F}  \ra \mc $ be a functor of $\mm$-modules. Then $F$ is a left adjoint if and only if it is
cocontinuous, in which case its right adjoint $G$ is
given by $$G(Y)_{k}=\mc (F(h_{k}),Y), \quad k\in K$$
\end{theorem}

\begin{example}\label{ex1}
Let $L$ be a small category and $K \st \Phi \ra L$ a $\mm$-functor. $\Phi$ induces a cocontinuous functor of $\mm$-modules.  $\mm^L \st F \ra \mm^K$. By Theorem \ref{nt1}, $F$ has a right adjoint functor $G$ given by
$$G(Y)_l=\mm^{K}(h_l\circ \Phi,Y).$$ 
\end{example}
\section{The closed $\mm$-model category structure on $\mc^K$}\label{s4}

We begin this Section with a brief review of the Reedy model category structure. Good references of the concept are \cite[Chapter III]{DHK}, \cite[Chapter 15]{HB}, \cite[Sections 5.1,  5.2]{H} and \cite{RV}. We also give the definitions of  closed symmetric monoidal model categories and model categories over them. We then study certain compatibility properties of Reedy diagrams in a model category over a monoidal model category.   \\

 A Reedy category is a small ordinary category $R$ equipped with two subcategories $R^{+}$ and $R^{-}$ and a function (called degree) $d:ob(R)\to\alpha$ where $\alpha$ is an ordinal number such that
\begin{enumerate}
	\item Every nonidentity morphism in $R^{+}$ (resp., $R^{-}$) raises (resp., lowers) degree.
	\item Every morphism $f$ in $R$  factors uniquely as $f=f^{+}\circ f^{-}$ where $f^{+}$ is a morphism in $R^{+}$
	and $f^{-}$ is a morphism in $R^{-}$.
\end{enumerate}
Observe that the subcategories $R^{-}$ and $R^{+}$ of $R$ are wide.
 $R$ is said to be direct (resp., inverse) if its subcategory  $R^{-}$ (resp., $R^{+}$) is the discrete category on the objects of $R$, in which case $R=R^{+}$ (resp., $R=R^{-})$. \\

Assume next that $R$ is a Reedy category and $\md$ is a model category. Let $\md^R$ be the category of functors from $R$ to $\md$.

 For $r\in R$, let $\delta(R^{+}/r)$   be the full subcategory of the over category $R^{+}/r$  whose objects are the nonidentity morphisms $s\to r$. 
Define $p^r:\delta(R^{+}/r)\to R$  to be the functor given by $p^r(s\to r)=s$.    For $X\in \md^R$, define $L_rX= \mbox{colim } X\circ p^r$. 

Similarly,  let $\partial(r/R^{-})$ be the full subcategory of the under category $r/R^{-}$ whose objects are the nonidentity morphisms  $r\to s$.
Define $p_r:\partial(r/R^{-})\to R$ to be the functor given by $p_r(r\to s)=s$.  For $X\in \md^R$, define $M_rX = \mbox{lim }  X\circ p_r$.\\

 The main theorem about Reedy categories says that there is a model category structure on the diagram category $\md^R$ with

\begin{itemize}
\item Weak equivalences are the objectwise weak equivalences.
\item Cofibrations are maps $A\to B$ such that the maps $A_r\coprod_{L_rA}L_rB\to B_r$ are cofibrations for all $r$.
\item Fibrations are maps $X\to Y$ such that the maps $X_r\to M_rX\times_{M_rY}Y_r$  are fibrations for all $r$.
\end{itemize}
The inclusion functor $R^+\to R$ induces a restriction functor   $T_{\md}:\md^R\to\md^{R^+}$ which preserves and reflects cofibrations and weak equivalences. By Example \ref{ex1},  $T_{\md}$ is left adjoint. It follows that $T_{\md}$ is a left Quillen functor.\\

We next recall the definitions of a symmetric monoidal model category \cite[Definition 4.2.6.]{H} and a $\mm$-model category \cite[Definition 4.2.18.]{H}. 
\begin{defn}
A closed symmetric monoidal model category is a closed symmetric monoidal category $\mm$ with monoidal product $\otimes$ which is also a model category such that the following axioms hold.
\begin{enumerate}
\item (Unit axiom) For any cofibrant approximation $Q1 \st{q} \ra 1$ of the unit $1$ of $\mm$ and any cofibrant object $m\in \mm$, the map  $Q1\otimes m \st{q\otimes m} \ra m$ is a weak equivalence.
\item (Pushout-product axiom) 
Let $f:a \ra b$ and $g:x \ra y$ be cofibrations in  $\mm$,  then the induced map $$f  \square  g : a \ot y \coprod_{a\ot x} b \ot x \ra b\ot y$$ is a cofibration in $\mm$ which is trivial if $f$ or $g$ is.
\end{enumerate}
\end{defn}
\begin{remark}\label{r3}
Let $(\mm,\otimes)$ be a closed symmetric monoidal category which is a model category .
\begin{enumerate}
\item By \cite[Lemma 4.2.2.]{H}, $\mm$ satisfies the pushout-product axiom  if and only if for every cofibration $f:a\ra b$ and every fibration $g:m\longrightarrow n$ in $\mm$, the induced map 
$$m^{b}\ra m^{a} \times_{n^{a}} n^{b}$$  
is a fibration in $\mm$ which is trivial if $f$ or $g$ is.
\item Assume that $\mm$ is cofibrantly generated. Then by \cite[Corollary 4.2.5.]{H}, it can easily be seen that
 $\mm$ satisfies the pushout-product axiom  if and only if for every generating (trivial) cofibration $f:a\ra b$ and every fibration $g:m\longrightarrow n$ in $\mm$, the induced map 
$$m^{b}\ra m^{a} \times_{n^{a}} n^{b}$$  
is a fibration in $\mm$ which is trivial if $f$ or $g$ is.\end{enumerate} 
\end{remark}  

\begin{defn}
Let $\mm$ be a closed symmetric monoidal model category. A 
$\mm$-model category is a closed $\mm$-module $\mc$ which is also a model category such that the following axioms hold 
\begin{enumerate}
\item (Unit axiom) For any cofibrant approximation $Q1 \st{q} \ra 1$ of the unit $1$ of $\mm$ and any cofibrant object $c\in \mc$, the map  $Q1\otimes c \st{q\otimes c} \ra c$ is a weak equivalence.
 \item (Pushout-product axiom) Let $f:a \ra b$ be a cofibration in  $\mm$ and $g:x \ra y$ a cofibration in  $\mc$,  then the induced map $$f  \square  g : a \ot y \coprod_{a\ot x} b \ot x \ra b\ot y$$ is a cofibration in $\mc$ which is trivial if $f$ or $g$ is.
\end{enumerate}  
\end{defn}
\begin{remark}\label{r4}
Let $(\mm,\otimes)$ be a closed symmetric monoidal model category and $\mc$ a closed $\mm$-module which is a model category.
\begin{enumerate}
\item By \cite[Lemma 4.2.2.]{H}, the action of $\mm$ on $\mc$ satisfies the pushout-product axiom  if and only if for every cofibration $f:a\ra b$ in $\mm$ and every fibration $g:x\longrightarrow y$ in $\mc$, the induced map 
$$y^{b}\ra y^{a} \times_{x^{a}} x^{b}$$  
is a fibration in $\mc$ which is trivial if $f$ or $g$ is.
\item By \cite[Lemma 4.2.2.]{H}, the action of $\mm$ on $\mc$ satisfies the pushout-product axiom  if and only if for every cofibration $f:a\ra b$ in $\mc$ and every fibration $g:x\longrightarrow y$ in $\mc$, the induced map 
$$\mc(b,y) \ra \mc(a,y) \times_{\mc(a,x)} \mc(b,x)$$  
is a fibration in $\mm$ which is trivial if $f$ or $g$ is.
\item Assume that $\mm$ is cofibrantly generated. Then by \cite[Corollary 4.2.5.]{H}, it can easily be seen that $\mc$ is a $\mm$-model category if and only if for every generating (trivial) cofibration $f:a\ra b$ in $\mm$ and every fibration $g:x\longrightarrow y$ in $\mc$, the induced map 
$$y^{b}\ra y^{a} \times_{x^{a}} x^{b}$$  
is a fibration in $\mc$ which is trivial if $f$ or $g$ is.
\end{enumerate} 
\end{remark}  

Let $K$, $\mm$ and $\mc$ be as in the previous Section. An argument similar to the one used at the beginning of Section \ref{s2} to prove that $\mm^K$ is a $\mm$-module can be used here to show that the following functors 
$$\begin{array}{rrcllll}
 \otimes:&\mm\ti \mc^K & \to & \mc^K& & &\\
 &(m, C) & \longmapsto  & m\otimes C:&K&\to&\mc\\
  &      &              &             &k&\longmapsto& m\otimes C_k\\
																									
\end{array}$$
 
$$\begin{array}{rrcllll}
 (-)^{(-)}:&\mm^{op}\ti \mc^K & \to & \mc^K& & &\\

 &(m, C) & \longmapsto  & C^m:&K&\to&\mc\\
  &     &              &             &k&\longmapsto& C_k^m\\
																					
\end{array}$$
$$\begin{array}{rrcl}
 \mc^K(-,-):&(\mc^K)^{op}\times \mc^K & \to & \mm\\
 &(C, D) & \longmapsto  &\mc^K(C,D)=\int_{k\in K}\mc(C_k,D_k)
 \end{array}$$
define a closed $\mm$-module structure on $\mc^K$.\\
 
Assume now that $K$ is a Reedy category, $\mm$ is a closed  symmetric monoidal model category and $\mc$ is a $\mm$-model category. Let $\mm ^K$ and $\mc ^K$ have the Reedy model category structure. Then one has the following two results, the first of which is due to Barwick  \cite[Corollary 3.37]{Ba}. We here give a slightly different proof.

\begin{proposition}\label{p1}
$\mc^K$ is a $\mm$-model category.
\end{proposition}
\begin{proof} 
Step 1. $K$ is a direct category.\\
To prove the unit axiom, observe that by  \cite[Remark 5.1.7]{H}, if an object $Z$ of $\mc^{K}$ is cofibrant, then the objects $Z_k$ are cofibrant in $\mc$. The unit axiom for the action of $\mm$ on $\mc^{K}$ then follows from the unit axiom for the action of $\mm$ on $\mc$. We next prove the Pushout-product axiom.\\
Let $a \st f \ra b $ be a cofibration in $\mm$ and $X  \st g \ra Y$   a fibration in $\mc^{K}$, then the maps $X_{k} \st {g_{k}}
\ra  Y_{k}$ are fibrations in $\mc$ which are trivial if $g$ is.
Therefore the maps $$X_{k}^{b} \ra X_{k}^{a} \times _{Y^{a}_{k}} Y_{k}^{b}$$
are fibrations in $\mc$ which are trivial if $ f $  or $g$ is. It follows that
$X^{b} \ra X^{a} \times _{Y^{a}} Y^{b}$ is a fibration in
$\mc^{K}$ which is trivial if $f$ or $g$ is. The result then follows from the first point in Remark \ref{r4}. \\ 
Step 2. The general case.\\
Let $K^+$ be the direct subcategory of $K$ used to define the Reedy category structure on $K$. The inclusion functor $K^+\to K$ induces the restriction functor   $T_{\mc}:\mc^K\to\mc^{K^+}$ which  preserves and reflects weak equivalences and (trivial) cofibrations.\\
Let $Q1  \st u \ra 1$ be a cofibrant replacement for the unit $1$ of $\mm$ and $C$ a cofibrant object in $\mc^{K}$, $T_{\mc}(C)$ is cofibrant in $\mc^{K^{+}}$. By step 1,  $u\otimes T_{\mc}(C)$ is a weak equivalence, observe that $u\otimes T_{\mc}(C)=T_{\mc}(u\otimes C)$ and $T_{\mc}$ reflects weak equivalences, therefore $u\otimes C$ is a weak equivalence and the unit axiom holds.\\
For $f:a \ra b$ in $\mm$   and $g:X \ra Y$ in $\mc^K$ (resp., $\mc^{K^+}$), let $f   \square  g$ be the pushout-product
 $$f  \square  g : a \ot Y \coprod_{a\ot X} b \ot X \ra b\ot
Y$$ 
Assume that $f$  is a cofibration in $\mm$   and $g$ is a cofibration in $\mc^K$. $T_{\mc}(g)$ is a cofibration which is trivial if $g$ is. 
 By Step 1, $f\square T_{\mc}(g)$ is a cofibration which is trivial if $f$ or $g$ is. Observe that
$$f \square T_{\mc}(g)=T_{\mc}(f \square g)$$  Thus $f \square g$ is a cofibration which is trivial  if $f$ or $g$ is.
\end{proof}
\begin{example}
If $\mc$ is a simplicial model category and $K$ is a Reedy category. Then with the Reedy model structure, $\mc^K$ is a simplicial model category. 
\end{example}
The following Lemma was proved by Riehl and Verity in the case where $\mm$ is the category of simplicial sets \cite[Theorem 10.3.] {RV}.
\begin{lemma}\label{l7} 
The adjunction   $(\otimes,hom_l,hom_r)$ given by Lemma \ref{l3} is a Quillen adjunction of two variables from $\mm^K\times \mc$ to $ \mc^K$.
\end{lemma}
\begin{proof} 
Step 1: $K$ is direct. \\
	Let $f:a\to b$ be a cofibration in $\mc$ and $g:X\to Y$ be a fibration in $\mc^K$. $\mc$ is a closed $\mm$-model category and $g_k:X_k\to Y_k$ is a fibration in $\mc$, thus by \cite[Lemma 4.2.2.]{H}, the map 
	$$\mc(b,X_k)\to \mc(a,X_k)\times_{ \mc(a,Y_k)}\mc(b,Y_k)$$ is a fibration which is trivial if $f$ or $g$ is. Therefore the map
	$$hom_r(b,X)\to hom_r(a,X)\times_{ hom_r(a,Y)}hom_r(b,Y)$$
is a fibration in $\mc^K$ which is trivial if $f$ or $g$ is. The result now follows from \cite[Lemma 4.2.2.]{H}. \\ 
Step 2: The general case.\\
	Let $K^+$ be the direct subcategory of $K$ used to define the Reedy category structure on $K$. The inclusion functor $K^+\to K$ induces  restriction functors $T_{\mm}:\mm^K\to\mm^{K^+}$ and $T_{\mc}:\mc^K\to\mc^{K^+}$ which preserve and reflect (trivial) cofibrations.\\
For $f:A\to B$ in $\mm^K$ (resp., $\mm^{K^+}$) and $g:c\to d$ in $\mc$, let $$f\square g:A\otimes d\coprod_{A\otimes c} B\otimes c\to B\otimes d$$ be the the pushout-product of $f$ and $g$. Suppose that $f $ is a cofibration in $\mm^K$ and $g$ is a cofibration in $\mc$.
	$$T_{\mc}(f\square g)=T_{\mm}(f)\square g$$
	$T_{\mm}$ preserves (trivial) cofibration, by Step 1, $T_{\mm}(f)\square g$ is a cofibration which is trivial if $f$ or $g$ is. $T_{\mc}$ reflects (trivial) cofibration, thus $f\square g$ is a cofibration which is trivial if $f$ or $g$ is.
It follows that the adjunction $(\otimes,hom_l,hom_r)$ is a Quillen adjunction.
\end{proof}
\section{The monoidal model structure on $\mm^K$}\label{ns4}
In this Section, we prove that the category of functors from a small category to a bicomplete closed symmetric monoidal category is still a closed symmetric monoidal category. We then similarly prove that the category of functors from a Reedy category into a closed symmetric monoidal model category is again a closed symmetric monoidal model category, provided that either the unit is cofibrant or the ground category is cofibrantly generated.\\

Still in this Section, $K$ is assumed to be a small category and  $\mm$ a bicomplete closed symmetric monoidal category with monoidal product $\ot$ and unit $1$.
The functor 
$$\begin{array}{rrcllll}
 \otimes:&\mm^K\times \mm^K & \to & \mm^K& & &\\
         &(M, N) & \longmapsto  & M\otimes N:&K&\to&\mm\\
         &      &              &             &k&\longmapsto& M_k\otimes N_k\\
	\end{array}$$
defines a symmmetric monoidal structure on $\mm^K$.\\

The next result is due to Day \cite[Section 5]{D}. 
\begin{proposition}\label{nc3} 
$(\mm^{K}, \ot)$ is a closed symmetric monoidal category with internal homomorphism functor 
$$\begin{array}{ccllll}
(\mm^K)^{op}\times \mm^K & \to & \mm^K& & &\\
        (N, P) & \longmapsto  & P^N:&K&\to&\mm\\
              &              &             &k&\longmapsto& \mm^K (h_{k} \ot N,P)\\
\end{array}$$
\end{proposition}
\begin{proof}
For $N\in \mm^{K}$ define a functor $$F_{N}:\mm^{K} \ra \mm ^{K}
$$
                                                       $$M\longmapsto M \ot N$$
$F_{N} $ is a cocontinuous functor of $\mm$-modules, by Theorem \ref{nt1}, $F_{N}$ has a left adjoint functor $G_N$ given by
$G_N(P)_{k}=\mm^K (h_{k} \ot N,P)$.
\end{proof}
\begin{lemma}\label{nl5}
Assume that $K$ is a direct category, $\mm$ a closed symmetric monoidal model category and let $\mm^K$ have the Reedy model category structure. 
\begin{enumerate}
\item Let $p, q \in K$ and $a \st f\ra b$ be a (trivial) cofibration in $\mm$, then the map 
$$h_{p} \ot h_{q} \ot f: h_{p} \ot h_{q} \ot a  \ra h_{p} \ot h_{q} \ot b $$ 
is a (trivial) cofibration in $\mm^K$. 
\item Assume that the unit $1$ of $\mm$ is cofibrant and let 
$q \in K$ and  $A \st F\ra B$ be a (trivial) cofibration in $\mm^K$, then the map 
$$ h_{q} \ot F: h_{q} \ot A  \ra   h_{q} \ot B $$ 
is a (trivial) cofibration in $\mm^K$. 
\end{enumerate} 
\end{lemma}

\begin{proof} 
\noindent
\begin{enumerate}
 \item  For $k \in K$, let $K_0((p,q),k)$ be the subset of $K_0(p,k)\times K_0(q,k)$ of elements $(\alpha, \beta)$ such that there exists  

\[
   \left\{
                \begin{array}{l}
                  l \in K ,  l \neq k \\
                  \theta \in K_0(l,k)	\\
                  (\alpha^{'} ,  \beta^{'})\in K_0(p,l)\times K_0(q,l)
                \end{array}
              \right.
  \]
with
\[
   \left\{
                \begin{array}{l}
                  \alpha =\theta o \alpha^{'} \\
                  \beta =\theta o \beta^{'}
                \end{array}
              \right.
  \]
An easy inspection shows that the latching object $L_k (h_{p}\ot h_{q}) $ is given by 

 $$L_{k}(h_{p}\ot h_{q})\cong \coprod_ {K_0((p,q),k)} 1$$ and $$h_{p}(k)\ot h_{q}(k) \cong \coprod_ {K_0(p,k)\times I(q,k)} 1$$  where $1$ is the unit of $\mm$. Observe that 
$$L_{k}(h_{p}\ot h_{q}\ot a)\cong L_{k}(h_{p}\ot h_{q})\ot a$$ 
Thus the map 
$$(h_{p}\ot h_{q}\ot a)_k\coprod_ {L_{k}(h_{p}\ot h_{q}\ot a)}L_{k}(h_{p}\ot h_{q}\ot b)\ra (h_{p}\ot h_{q}\ot b)_k$$ is simply the map
$$(\coprod_{K_0((p,q),k)} b)\coprod (\coprod_{K_0(p,k)\times K_0(q,k)-K_0((p,q),k)} a)\ra (\coprod_{K_0((p,q),k)} b)\coprod (\coprod_{K_0(p,k)\times K_0(q,k)-K_0((p,q),k)}b)$$
which is the identity on the first summand and the iteration of $f$ on the second summand  and is clearly a (trivial) cofibration.
\item $\mm^{K}$ is closed, therefore the property that $ h_{q} \ot F: h_{q} \ot A  \ra   h_{q} \ot B $ is a (trivial) cofibration  whenever
$F:  A \ra B $ is a (trivial) cofibration is equivalent to the property that  $g^{h_{q}}:  M^{h_{q}} \ra N^{h_{q}}$ ia a (trivial) fibration whenever $g:  M \ra N$ is a (trivial)  fibration. So let $g:  M \ra N $ be a (trivial) fibration and let $p\in K$. The unit $1$ of $\mm$ cofibrant, therefore by the first part of this Lemma, $h_{p}\ot h_{q}$ is cofibrant. By Proposition \ref{p1}, $\mm^K$ is a $\mm$-model category, thus
$$\mm^K(h_{p}\ot h_{q},M)\longrightarrow \mm^K(h_{p}\ot h_{q},N)$$
ia a (trivial) fibration. This means that $$g^{h_{q}}:  M^{h_{q}} \ra N^{h_{q}}$$ is an objectwise (trivial) fibration (Proposition \ref{nc3}). $K$ is direct, therefore $g^{h_{q}}$ is a (trivial) fibration.  
\end{enumerate} 
\end{proof}

\begin{remark}\label{r2} Assume that $K$ is a direct category and $\mm$  is a cofibrantly generated closed symmetric monoidal model category with monoidal product $\otimes$. By Lemma \ref{nl3}, the q-evaluation functor $\mm^{K} \st {Ev_{q}}\ra\mm$ is right adjoint to the functor $\mm \st {F_{q}}\ra \mm^{K}$ given by $F_{q}(m)=h_{q}\ot m$. Therefore by \cite[Remark 5.1.8.]{H}, $\mm^{K}$ is a cofibrantly generated model category with generating (trivial) cofibrations are the maps 
$$F_{q}(f)=h_{q}\ot f:h_{q}\ot a \ra h_{q}\ot b $$ 
where $f:a \ra b$ is a generating (trivial) cofibration in $\mm$.
\end{remark}

We are now ready to prove our first main result of which Barwick's Theorem \cite[Theorem 3.51]{Ba} is a special case.

\begin{theorem}\label{nt2}
Let $K$ be a Reedy category and $\mm$ a closed symmetric monoidal model category. Assume further that either the unit $1$ of $\mm$ is cofibrant or $\mm$ is cofibrantly generated. Then with the Reedy model category structure, $\mm^{K}$ is a closed symmetric monoidal model category.
\end{theorem}
\begin{proof}
By Proposition \ref{nc3} $\mm^{K}$ is a closed symmetric monoidal category. The unit axiom for $\mm^{K}$ follows easily from an argument similar to the one used in   the proof of Proposition \ref{p1}. It remains to prove the Pushout-product axiom. \\
Step 1: $K$ is direct.\\
Case 1: The unit $1$ is cofibrant.\\
Let $f:A\ra B$ be a cofibration in $\mm^K$, $g:M\rightarrow N$  be a fibration $\mm^K$ and let $q \in K$. By the second point of Lemma \ref{nl5}, the map
$$ h_{q} \ot f: h_{q} \ot A  \ra   h_{q} \ot B $$
is a cofibration $\mm^K$ which is trivial if $f$ is. By Proposition \ref{p1}, $\mm^K$ is a $\mm$-model category, therefore by the second point of Remark \ref{r4}, the map 
 
 $$\mathcal{V}^{K}(h_{q}\otimes B,M)\longrightarrow \mathcal{V}^{K}(h_{q}\otimes A,M)\times_{\mathcal{V}^{K}(h_{q}\otimes A,N)}\mathcal{V}^{K}(h_{q}\otimes B,N)$$ 
is a fibration in $\mm$ which is trivial if $f$ or $g$ is. $K$ is direct, thus  the map  $$M^{B}\longrightarrow M^{A}\times_{N^{A}}N^{B}$$
is a  fibration in $\mathcal{V}^{K}$  which is trivial if $f$ or $g$ is.\\ 
Case 2: $\mm$ is cofibrantly generated.\\
Let $f:a\ra b$ be a generating (trivial) cofibration in $\mm$, $g:M\rightarrow N$  be a fibration  in $\mm^K$ and $p,q\in K$. By the first point of Lemma \ref{nl5}, the map
$$ h_{p} \ot h_{q} \ot f: h_{p} \ot h_{q} \ot a  \ra h_{p} \ot h_{q} \ot b $$ is a cofibration which is trivial if $f$ is.
By Proposition \ref{p1}, $\mm^K$ is a $\mm$-model category, therefore by the second point of Remark \ref{r4}, the map 
$$\mm^{K}(h_{p}\ot h_{q}\ot b,M)\rightarrow \mm^{K} (h_{p}\ot h_{q} \ot a,M)\times_{\mm^{K}(h_{p} \ot h_{q} \ot a,N)}\mm^{K}(h_{p}\ot h_{q}\ot b,N)$$ is a fibration in $\mm$ which is trivial if $f$ or $g$ is. By Proposition \ref{nc3}, the map
$$M^{h_{q} \ot b}\ra M^{h_{q} \ot a} \times_{N^{h_{q} \ot a}}
N^{h_{q} \ot b}$$ 
 is an objectwise fibration  which is trivial if $f$ or $g$ is. $K$ is direct, thus the map $M^{h_{q} \ot b}\ra M^{h_{q} \ot a} \times_{N^{h_{q} \ot a}}
N^{h_{q} \ot b}$ is  a fibration which is trivial if $f$ or $g$ is. By Remark \ref{r2} and the second point in Remark \ref{r3}, $\mm^{K}$ is a  symmetric monoidal model
category.\\
Step 2: The general case. \\		
The result is an immediate consequence of step 1 and the fact that the restriction functor $T_{\mm}:\mm^K \ra \mm^{K^{+}}$ is monoidal, preserves pushouts and preserves and reflects (trivial) cofibrations.
 \end{proof}
\section{The closed $\mm^K$-model structure on $\mc^K$}\label{s5}
Let $K$, $\mm$ and $\mc$ be as in Section \ref{s2}. We here prove our final main results. More precisely, we show that $\mc^K$ is a closed $\mm^K$-module. We then further show that when $K$ is a Reedy category, $\mm$ a closed symmetric monoidal model category which either has a cofibrant unit or is cofibrantly generated and $\mc$ a closed $\mm$-model category, then with the Reedy model category structures, $\mc^{K}$ is a closed $\mm^{K}$-model category.\\

Let $hom_l$ the functor defined in Section \ref{s2} and $\mc^K(-,-)$ the functor defined in Section \ref{s4}.  Define
$$\begin{array}{rrcllll}
 \otimes:&\mm^K\times \mc^K & \to & \mc^K& & &\\
         &(M, C) & \longmapsto  & M\otimes C:&K&\to&\mc\\
         &      &              &             &k&\longmapsto& M_k\otimes C_k\\
	\end{array}$$
$$\begin{array}{rrcllll}				
	(-)^{(-)}:&(\mm^K)^{op}\times\mc^K&\to&\mc^K&&&\\
	           &(M,X)&\longmapsto&X^{M}:&K&\to&\mc\\
	           &      &              &          &k&\longmapsto&\mbox{hom}_l(M\otimes h_k,X)\\
\end{array}$$
$$\begin{array}{rrcllll}		
	\mbox{Hom}:&(\mc^K)^{op}\times\mc^K&\to&\mm^K&&&\\
		          &(X,Y)&\longmapsto&\mbox{Hom}(X,Y):&K&\to&\mm\\
							&      &           &                        &k&\longmapsto&\mc^K(h_k\otimes X,Y)\\
\end{array}$$
\begin{theorem}\label{t0}
The above three functors define a closed $\mm^K$-module structure on $\mc^K$.
\end{theorem}
\begin{proof}
Let $M\in\mm^K$, $X$, $Y\in\mc^K$
$$\begin{array}{rcll}
\mm^K_0(M,\mbox{Hom}(X,Y))&\cong&\int_{k\in K}\mm_0(M_k,\mbox{Hom}(X,Y)_k)&\\
                                 &\cong&\int_{k\in K}\mm_0(M_k,\mc^K(h_k\otimes X,Y))&\\
																&\cong&\int_{k\in K}\mc^K_0( M_k\otimes h_k\otimes X,Y) & \\
                                  &\cong&\int_{k\in K}\mc^K_0(  h_k\otimes M_k\otimes X,Y)&\\
																	&\cong&\mc^K_0(\int^{k\in K} h_k\otimes M_k\otimes X,Y)&\\
																	&\cong&\mc^K_0((\int^{k\in K} h_k\otimes M_k)\otimes X,Y)&\\
																	&\cong&\mc^K_0(M\otimes X,Y)&\mbox{(by Yoneda coreduction \ref{l0})} \\
		&&&\\															
\mc^K_0(X,Y^{M})&\cong&\int_{k\in K}\mc_0(X_k,(Y^{M})_k)&\\
                         &\cong&\int_{k\in K}\mc_0(X_k,\mbox{hom}_l(M\otimes h_k,Y))&\\
												&\cong&\int_{k\in K}\mc^K_0(M\otimes h_k\otimes X_k,Y)& \mbox{(by Lemma \ref{l3})}\\
										&\cong&\mc^K_0(\int^{k\in K} M\otimes h_k\otimes X_k,Y)&\\	
												&\cong&\mc^K_0(M\otimes\int^{k\in K}h_k\otimes X_k,Y)&\\	
								&\cong&\mc^K_0(M\otimes X,Y)&\mbox{(by Yoneda coreduction \ref{l5})}
\end{array}	$$
\end{proof}
\begin{example}
Let $\Delta$ be the simplex category, $\ms$ the category of sets and assume that $\mc$ is a bicomplete category. $\ms$ with the cartesian product is a closed symmetric monoidal category and $\mc$, being bicomplete, is a closed $\ms$-module. Therefore,
by Theorem \ref{t0}, the category $s\mc$ of simplicial objects in $\mc$ is a closed simplicial category while the category $\mc^{\Delta}$ of cosimplicial objects in $\mc$ is a closed module over the category of cosimplicial sets.
\end{example}
We are now ready to prove the second main result of the paper.

\begin{theorem}\label{t1}\hspace{1mm}
Let $K$ be a Reedy category, $\mm$ a  closed symmetric monoidal model category and $\mc$ a closed $\mm$-model category. Assume further that either the unit $1$ of $\mm$ is cofibrant or $\mm$ is cofibrantly generated. Then with the Reedy model category structures, $\mc^{K}$ is a closed $\mm^{K}$- model category.
\end{theorem}
\begin{proof}
By Theorem \ref{t0}, $\mc^{K}$ is a closed $\mm^{K}$-module. By \cite[Remark 5.1.7]{H} and its dual, if $X$ is a cofibrant object in $\mc^{K}$, then the objects  $X_k$ are cofibrant in $\mc$ and a trivial fibration in $\mm^{K}$ is an objectwise trivial fibration in $\mm$ (the converse is false). The unit axiom for the action of $\mm^{K}$ on $\mc^{K}$ then follows from the unit axiom for the action of $\mm$ on $\mc$. We next prove the Pushout-product axiom.\\ 
Step 1: $K$ is direct.\\
Case 1: The unit $1$ is cofibrant.\\
Let $f:A\ra B$ be a cofibration in $\mm^K$, $g:X\rightarrow Y$  be a fibration $\mc^K$ and let $q \in K$. $\mm^K$ is symmetric, thus by the second point of Lemma \ref{nl5}, the map
$$ f \ot h_{q}:   A \ot h_{q} \ra B \ot h_{q} $$
is a cofibration which is trivial if $f$ is. By Lemma \ref{l7} 
 $$hom_l(B \ot h_{q},X)\longrightarrow hom_l (A \ot h_{q},X)\times_{hom_l(A \ot h_{q},Y)}hom_l(B \ot h_{q},Y)$$ 
is a fibration which is trivial if $f$ or $g$ is. This means that the map  $$X^{B}\longrightarrow X^{A}\times_{Y^{A}}Y^{B}$$
is a fibration in $\mc^{K}$  which is trivial if $f$ or $g$ is. By the first point of Remark \ref{r4}, the action of $\mm^K$ on $\mc^K$satisfies the Pushout-product axiom.\\
Case 2: $\mm$ is cofibrantly generated.\\
Let $f:a\ra b$ be a generating (trivial) cofibration in $\mm$, $g:X\rightarrow Y$  be a fibration  in $\mc^K$ and $p,q\in K$.\\ By the first point of  Lemma \ref{nl5}
$$h_{q}  \ot f \ot h_{p}:h_{q}  \ot a \ot h_{p}  \ra h_{q}  \ot b \ot h_{p} $$ 
is a cofibration in $\mm^K$ which is trivial if $f$ is trivial. By Lemma \ref{l7}  
$$hom_l(h_{q}\ot b \ot h_{p},X)\ra hom_l (h_{q} \ot
a \ot h_{p},X)\times_{hom_l(h_{q} \ot a \ot h_{p},Y)}hom_l(h_{q}\ot b \ot h_{p},Y)$$
is a fibration in $\mc$ which is trivial if $f$ or $g$ is.  It follows that the map

$$ (X^{h_{q} \ot b})_p\ra (X^{ h_{q} \ot a })_p \times_{ (Y^{h_{q} \ot a})_p}(Y^{h_{q} \ot b})_p$$
is a fibration in $\mc$ which is trivial if $f$ or $g$ is. $K$ is direct, therefore the fibrations in $\mc^{K}$ are the objectwise fibrations. Hence
$$ X^{h_{q} \ot b}\ra X^{ h_{q} \ot a } \times_{ Y^{h_{q} \ot a}}Y^{h_{q} \ot b}$$
is a fibration in $\mc^K$ which is trivial if $f$ or
$g$ is. By Remark \ref{r2} and the third point of Remark \ref{r4},  the action of $\mm^K$ on $\mc^K$ satisfies the Pushout-product axiom.\\ 
Step 2: $K$ is an arbitrary Reedy category.\\
An argument that uses step 1 and is strictly similar to the one used in step 2 in the  proof of  Theorem \ref{nt2} shows the desired result.
\end{proof}
\section{Other properties and examples}\label{s6}
Let $\mc, \md$ be two closed $\mm$-modules and $\oppairi{\mc}{\md}{F}{G}$ an adjunction. By  \cite[Theorem 4.85]{K}, $F$ is a functor of $\mm$-modules if and only if $G^{op}:\md^{op}\to\mc^{op}$ is a functor of $\mm$-modules with the dual action of $\mm$ on $\mc^{op}$ and $\md^{op}$. If that is the case then the adjunction $(F,G)$ is called an adjunction of $\mm$-modules.
\begin{defn}
Let $\mm$ be a closed symmetric monoidal model category and $\mc, \md$ closed $\mm$-model categories. An adjunction
$\begin{tikzcd}
\mc\ar[r,shift left=.75ex,"F"]
  \ar[r,shift right=.75ex,leftarrow,swap,"G"]
&
\md
\end{tikzcd}$
is called a Quillen adjunction (resp., Quillen equivalence) of $\mm$-modules if
\begin{enumerate}
	\item $(F,G)$ is an adjunction of $\mm$-modules.
	\item $(F,G)$ is a Quillen adjunction (resp., Quillen equivalence).
\end{enumerate}
\end{defn}
\begin{theorem}\label{t2}
Let $\mm$ be a closed symmetric monoidal model category which either has a cofibrant unit or is cofibrantly generated, $\begin{tikzcd}
\mc\ar[r,shift left=.75ex,"F"]
  \ar[r,shift right=.75ex,leftarrow,swap,"G"]
&
\md
\end{tikzcd}$ 
a Quillen adjunction (resp., Quillen equivalence) of $\mm$-model categories and $K$ a Reedy category. Then with the Reedy model category structures, the induced adjunction $\begin{tikzcd}
\mc^K\ar[r,shift left=.75ex,"F^K"]
  \ar[r,shift right=.75ex,leftarrow,swap,"G^K"]
&
\md^K
\end{tikzcd}$
is a Quillen adjunction (resp., Quillen equivalence) of $\mm^K$-model categories.
\end{theorem}
\begin{proof} Let $f:X\ra Y$ be a cofibration in $\mc^K$. Being left adjoint, $F$ is cocontinuous. It follows that the map
$$F^{K}(X)_k\coprod _{L_kF^{K}(X)}L_kF^{K}(Y)\to F^{K}(Y)_k$$ is just the map $$F(X_k\coprod _{L_k X}L_k Y\to Y_k).$$
$F$ preserves cofibrations. It follows that $F^{K}(f)$ is a cofibration and $F^{K}$ preserves cofibrations. The dual argument 
shows that $G^{K}$ preserves fibrations, therefore $(F^{K},G^{K})$ is a Quillen adjunction. Assume now that $(F,G)$ is a Quillen equivalence and let $A$ be a cofibrant object in $\mc^K$, $B$ be a fibrant object in $\md^K$ and $F^K(A) \st{g}\ra B$ a map in $\md^{K}$ with adjoint $ A \st{g^{\#}}\ra G^K(B)$. The  objects $A_k$ are  cofibrant in $\mc$ the  objects $B_k$ are fibrant in $\md$, therefore the maps  $F(A_k) \st{g_k}\ra B_k$ are weak equivalences if and only if the maps  $ A_k \st{g_k^{\#}}\ra G(B_k)$ are weak equivalences. i.e., $g$ is a weak equivalence if and only if $g^{\#}$ is a weak equivalence. It follows that $(F^K,G^K)$ is a Quillen equivalence.
\end{proof}
\begin{example}
Let $\mm$ be a bicomplete closed symmetric monoidal category with monoidal product $\otimes$, unit $1$, and internal hom functor denoted exponentially so that there is a natural bijection
$$\mm_0(x\otimes y,z)\cong\mm_0(x,z^y) \qquad \mbox { for } x,y,z\in\mm$$
 Let $a$ be monoid a in $\mm$ (see \cite[VII.3]{ML}) and $\mm_a$ the category of right $a$-modules in $\mm$. An object in $\mm_a$ consists of an object $m\in \mm$ together with an action map  
$m\otimes a \st {\lambda_m} \ra m$
satisfying the associativity and the unit conditions. A map in $\mm_a$ is an arrow $m\ra n$ in $\mm$ which respects the actions of $a$ on $m$ and $n$. Precise definitions may be found in \cite [VII.4]{ML}. The monoidal product in $\mm$ induces a $\mm$-module structure on $\mm_a$
$$ \begin{array}{rrl}
 \otimes:&\mm\times\mm_a \to & \mm_a\\
 &(x,m)\longmapsto  & x\otimes m																
\end{array}$$
For $m,n\in\mm_a$, the set of morphisms $(\mm_a)_0(m,n)$ is the equalizer
\begin{equation}\label{eq4}
(\mm_a)_0(m,n)\to\begin{tikzcd}
(\mm)_0(m,n)\ar[r,shift left=.75ex,""]
  \ar[r,shift right=.75ex,rightarrow,swap,""]
&
(\mm)_0(m\otimes a,n)
\end{tikzcd}
\end{equation}
where one of the parallel arrows is induced by the action of $a$ on $m$ and the other takes an arrow $m \st {f}\ra n$  to the composite $m\otimes a \st {f\otimes a} \ra n\otimes a\st {\lambda_n} \ra n$.\\
For $x,y\in \mm$, let $y^{x}\otimes x\st {\theta^{x}_{y}}\ra y$ be the adjoint of the identity map $y^{x}\ra y^{x}.$\\
Let $\mm_a(m,n)$ be the equalizer 
$$\mm_a(m,n)\to\begin{tikzcd}
n^m\ar[r,shift left=.75ex,""]
  \ar[r,shift right=.75ex,rightarrow,swap,""]
&
n^{m\otimes a}
\end{tikzcd}$$
where one of the parallel arrows is the adjoint of the composite
$$n^m\otimes (m\otimes a)\ra (n^m\otimes m)\otimes a \xrightarrow{{\theta^{m}_{n}\otimes a}} n\otimes a \st{\lambda_n} \ra n$$
and the other is the adjoint of the composite
$$n^m\otimes (m\otimes a)\xrightarrow{{n^m\otimes \lambda_m}} n^m\otimes m \st{\theta^{m}_{n}}\ra  n.$$
We then have a functor
$\mm_a(-,-):\mm_a^{op}\times\mm_a\to\mm$ satisfying $$(\mm_a)_0(x\otimes m,n) \cong \mm_0(x,\mm_a(m,n)) \qquad  \mbox{for }x\in \mm \mbox{  and } m,n \in\mm_a.$$
For $m\in\mm_a$, let $m\otimes a\to m$ be the action of $a$ on $m$, $m\to m^a$ its adjoint. For $x\in\mm$, the map $m\to m^a$ induces a map 
$$m^x\to (m^a)^x\cong m^{a\otimes x}\cong m^{x\otimes a}\cong (m^x)^a$$
its adjoint
$$m^x\otimes a\to m^x$$
induces an $a$-module structure on $m^x$, we therefore have a functor
$$ \begin{array}{rrl}
 \mm^{op}\times\mm_a &\to & \mm_a\\
 (x,m)&\longmapsto  & m^x															
\end{array}$$
An easy inspection shows that 
$$(\mm_a)_0(x\otimes m,n) \cong (\mm_a)_0(m,n^x)$$
$\mm_a$ is then a closed $\mm$-module. Assume that $\mm$ is a cofibrantly generated monoidal model category and the monoid $a$ is cofibrant in $\mm$. By \cite[Corollary 2.2.]{HH},
$\mm_a$ has a cofibrantly generated model category structure where a map is a fibration or a weak equivalence if it is so in $\mm$. Furthermore, any cofibration in $\mm_a$ is a cofibration in $\mm$.\\
	Let $x\to y$ be a cofibration in $\mm$ and $m\to n$ a fibration in $\mm_a$. The induced map $m^y\to m^x\times_{n^x} n^y$ is a fibration in $\mm_a$ which is trivial if $x\to y$ or $m\to n$ is. $\mm_a$ is then a $\mm$-model category. Let $b$ be another cofibrant monoid in $\mm$ and $a\st {f} \ra b$ a monoid map. For $m\in\mm_a$, let $m\otimes_ab$ be the coequalizer 
	$$
	\begin{tikzcd}
m\otimes a\otimes b\ar[r,shift left=.75ex,""]
  \ar[r,shift right=.75ex,rightarrow,swap,""]
&
{m\otimes b}
\end{tikzcd}
\to m\otimes_{a}b.$$
where the parallel arrows are induced by the right action of $a$ on $m$ and the left action of $a$ on $b$.
$m\otimes_a b$ is a right $b$-module and the functor
$$-\otimes_{a}b:\mm_a\to\mm_b$$
is a $\mm$-functor which is left adjoint to the restriction functor
$$Res:\mm_b\to\mm_a.$$
Since $Res$ preserves (trivial) fibrations, it follows that the adjunction $(-\otimes_ab,\Res)$ is a Quillen adjunction of $\mm$-model categories.  Furthermore, if $a\st {f} \ra b$ is a weak equivalence and the domains of the generating cofibrations in $\mm$ are cofibrant, then  by \cite[Theorem 2.4.]{HH}, $(-\otimes_ab,\Res)$ is a Quillen equivalence. It follows that the adjunction
$$\begin{tikzcd}
\mm_a^K\ar[r,shift left=.75ex,"-\otimes_ab^K"]
  \ar[r,shift right=.75ex,leftarrow,swap,"Res^K"]
&
\mm_b^K
\end{tikzcd}$$
is a Quillen adjunction of $\mm^K$-model categories. If further $f$ is a weak equivalence, then the above adjunction is a Quillen equivalence of $\mm^K$-model categories.
\end{example} 
\noindent{\bf Notation:}
For any model category $\ma$, let $s\ma =\ma^{\Delta^{op}}$ with the Reedy model category structure. Let $s^n\ma$ be the model category inductively defined by $s^0\ma=\ma$ and $s^{n+1}\ma=s(s^{n}\ma)$. Let $S$ be the model category of simplicial sets and define  $S^n=s^{n-1}S$, $n\geq 1$.
\begin{corollary}
Let $\begin{tikzcd}
\mc\ar[r,shift left=.75ex,"F"]
  \ar[r,shift right=.75ex,leftarrow,swap,"G"]
&
\md
\end{tikzcd}$
be a Quillen adjunction (resp., a Quillen equivalence) of simplicial model categories, then the induced adjunction
$\begin{tikzcd}
s^n\mc\ar[r,shift left=.75ex,"s^nF"]
  \ar[r,shift right=.75ex,leftarrow,swap,"s^nG"]
	&
s^n\md
\end{tikzcd}$
is a Quillen adjunction (resp., a Quillen equivalence) of $S^{n+1}$-model categories.
\end{corollary}
\begin{proof} $S$ is cofibrantly generated. Using induction and applying \cite[Proposition 7.7.]{RV}, we prove that $S^n$ is cofibrantly generated. The proof of the corollary is by induction and uses Theorem \ref{t2}. 
\end{proof}
\begin{example}
	 Let $\begin{tikzcd}
S\ar[r,shift left=.75ex,"|.|"]
  \ar[r,shift right=.75ex,leftarrow,swap,"Sing"]
&
Top
\end{tikzcd}$
be the usual realization-singular adjunction between $S$ and category $Top$ of compactly generated weakly Hausdorff topological spaces. $(|.|,\mbox{Sing})$ is a Quillen equivalence of simplicial model categories. Therefore $(s^n|.|,s^n\mbox{Sing})$ is a Quillen equivalence of $S^{n+1}$-model categories. 
\end{example}
\begin{example}
	Let $R$ be a commutative ring, $Mod_R$ the category of $R$-modules, $Alg_R$ the category of commutative $R$-algebras. The forget functor $O:Alg_R\to Mod_R$ has a left adjoint $T$ which is the symmetric algebra functor. By Proposition 4.2 and Theorem 4.17 in ~\cite{GS}, there are simplicial model category structures on $sMod_R$ and $sAlg_R$ so that 
 $$\begin{tikzcd}
sMod_R\ar[r,shift left=.75ex,"sT"]
  \ar[r,shift right=.75ex,leftarrow,swap,"sO"]
&
sAlg_R
\end{tikzcd}$$
is a Quillen adjunction of simplicial model categories, it follows that the adjunction 
$$\begin{tikzcd}
s^nMod_R\ar[r,shift left=.75ex,"s^nT"]
  \ar[r,shift right=.75ex,leftarrow,swap,"s^nO"]
&
s^nAlg_R
\end{tikzcd}$$
is a Quillen adjunction of $S^n$-model categories.
\end{example}
 \begin{remark}\label{r1} 
Though useful, Reedy categories are skeletal. They do not admit nonidentity isomorphisms and are not stable under equivalences. They do not therefore apply to many interesting cases occurring in topology. For these reasons, Berger and Moerdijk introduced the notion of generalized Reedy categories which overcome the above deficiencies, while keeping true the main theorem about model categories \cite{BM}. Examples of generalized Reedy categories which are generally not  strict Reedy categories include the category of  finite sets, the category of pointed finite sets, groupoids, Segal’s category, Connes’ cyclic category and others \cite[Examples 1.9]{BM}. We expect that the results of this paper extend to their categories.
\end{remark}

M. Ghazel\\
 Faculté des Sciences Mathématiques, Physiques et Naturelles de Tunis. University of Tunis El Manar, 2092 El Manar Tunis, Tunisia.\\
 moncef.ghazel@gmail.com\\
																	
\noindent F. Kadhi\\
Ecole Nationale des Sciences de l'Informatique, Manouba University. Manouba 2010, Tunisia.\\
					            fethi.kadhi@ensi-uma.tn
						
\end{document}